\documentclass{amsart}
\usepackage[leqno]{amsmath}

\usepackage{amsmath}
\usepackage{amsthm}
\usepackage{amssymb}
\usepackage{amscd}
\usepackage{amsxtra}     
\usepackage{epsfig}
\usepackage{verbatim}
\usepackage[all]{xypic}

\theoremstyle{plain} 

\newtheorem{thm}{Theorem}[section]
\newtheorem{lem}[thm]{Lemma}
\newtheorem{prop}[thm]{Proposition}

\newtheorem{cor}[thm]{Corollary}
\newtheorem{dfn}[thm]{Definition}
\newtheorem{eg}[thm]{Example}

\newtheorem{conj}[thm]{Conjecture}
\newtheorem{cons}[thm]{Consequence}

\theoremstyle{remark}
\newtheorem*{rmk}{Remark}

\numberwithin{equation}{section}

\newcommand{\PP}{\mathbf{P}}

\newcommand{\ZZ}{\mathbf{Z}}

\newcommand{\Div}{\textup{div}}

\newcommand{\tensor}{\otimes}

 \DeclareMathOperator{\Tor}{Tor}
 \DeclareMathOperator{\Ext}{Ext}
 \DeclareMathOperator{\Hom}{Hom}

 \DeclareMathOperator{\Supp}{Supp}
 \DeclareMathOperator{\Spec}{Spec}
 \DeclareMathOperator{\Sing}{Sing}
 
\DeclareMathOperator{\CH}{CH}
 \DeclareMathOperator{\Proj}{Proj}

 \DeclareMathOperator{\pd}{pd}

 \DeclareMathOperator{\depth}{depth}

 \DeclareMathOperator{\ch}{ch}
 \DeclareMathOperator{\syz}{syz}

 \newcommand{\Min}{\textup{Min}}
\newcommand{\ord}{\textup{ord}}

\newcommand{\Ann}{\textup{Ann}}

\begin{document}

\bibliographystyle{plain}

\title[Decent intersection and Tor-rigidity]{Decent intersection and Tor-rigidity for modules over local hypersurfaces}

\author{Hailong Dao}
\address{Department of Mathematics\\
University of Kansas\\
 Lawrence, KS 66045-7523 USA}
\email{hdao@math.ku.edu}

\thanks{The author is partially supported by NSF grant 0834050}


\maketitle
\begin{abstract}
We study two properties for a pair of finitely generated modules over a local hypersurface $R$:
decency, which is close to proper intersection of the supports,  and $\Tor$-rigidity. We show that the vanishing of Hochster's function $\theta^R(M,N)$, known to imply decent intersection, also
implies rigidity. We investigate the vanishing of $\theta^R(M,N)$
to obtain new results about decency and rigidity over
hypersurfaces. Many applications are given. 
\end{abstract}
\section{Introduction}

Throughout this paper we will deal exclusively with a local,
Noetherian, commutative ring $R$ and finitely generated modules
over $R$.

Two $R$-modules $M$ and $N$ such that $l(M\tensor_RN)<\infty $ are
said to \textit{intersect decently} if $\dim M + \dim N \leq \dim
R$. We say that $M$ is \textit {decent} if for all $N$ such that
$l(M\tensor_RN)<\infty $, $M$ and $N$ intersect decently. This
property arises naturally from Serre's work on intersection
multiplicity (\cite{Se}), which shows that over a regular local
ring, any two modules intersect decently. In fact, to have a
satisfying intersection theory over a local ring, one needs modules to
intersect decently as a minimum requirement. However, sufficient
conditions for decent intersection are much more elusive in
non-regular situation, even when $R$ is a hypersurface. For example,
the Direct Summand Conjecture would follow if one could show that a certain module over a local
hypersurface is decent (see \cite{Ho1}). As another example, it was conjectured by Peskine and Szpiro
(cf. \cite{PS}) that over any local ring, a module of finite projective dimension  is decent. This question is
open even for hypersurfaces of ramified regular local rings. Such questions serve as our initial motivation, and while our investigation did not settle any of those conjectures, we hope it will contribute towards better understanding of this intriguing condition. 

A pair of modules $(M,N)$ is called \textit{rigid} if for any integer $i\geq0$,
$\Tor_i^R(M,N)=0$ implies $\Tor_j^R(M,N)=0$ for all $j\geq i$. Moreover, $M$ is \textit{rigid} if for
all $N$, the pair $(M,N)$ is rigid. Auslander studied rigidity in order to understand torsion on
tensor products (\cite{Au}). He also observed that rigidity of $M$ implies other nice properties, such as
any $M$-sequence must be an $R$-sequence. To further demonstrate the usefulness of rigidity, let us recall the following depth formula (\cite{HW1}, 2.5):  suppose $R$ is a local complete intersection and $M,N$ are non-zero
finitely generated modules over $R$ such that $\Tor_i^R(M,N)=0$ for all $i\geq 1$. Then:
$$ \depth(M) + \depth(N) = \depth(R)+\depth(M\tensor_RN) $$

Thus, rigidity allows us to force a very strong condition on the
depths of the modules by proving the vanishing of one single
$\Tor$ module. Auslander's work, combined with results of Lichtenbaum
(\cite{Li}), showed that modules over regular local rings are
rigid. Huneke and Wiegand (\cite{HW1},\cite{HW2}) continued this
line to study rigidity over hypersurfaces, and their paper has major influence on 
our work here.


The classical condition conjectured to be sufficient for both
rigidity and decency was that one of the modules must have finite
projective dimension. In general, this is open for decency and
false for rigidity (see \cite{He}). In any case, having finite
projective dimension is too restrictive for the most interesting
applications, so a question arises: \textit{Are there more
flexible sufficient conditions for rigidity and decency?}  In this work, we obtain some answers to this question for modules over \textit{local hypersurfaces} using a key  function defined by Hochster (see Section \ref{2}), whose vanishing controls both
rigidity and decency, together with  some results from
Intersection Theory. Our results often
involve conditions about the classes of the modules in the
Grothendieck group of finitely generated modules over $R$ which are weaker than having finite
projective dimension.

Throughout this note we will assume that our hypersurface $R$
comes from an  unramified or equicharacteristic regular local ring
(we call such hypersurfaces ``admissible").
Since we need to apply results such as Serre's Positivity and
Non-negativity of $\chi_i$, which are open in general for the
ramified case, this is necessary. In some particular instances,
such as in low dimensions, this restriction can be relaxed,
however we feel it may disrupt the flow of the paper to comment on
every such case. The reader will lose very little by thinking of
the equicharacteristic (containing a field) case.

Section \ref{2} is a review of basic notation and some preliminary
results. Of particular interest is  Hochster's theta function. For
a local hypersurface $R$ and  a pair of finitely generated $R$
modules $M,N$ such that $l(\Tor_i^R(M,N))<\infty$ for all $i\gg0$,
we can define:
$$\theta^R(M,N) = l(\Tor_{2e+2}^R(M,N)) - l(\Tor_{2e+1}^R(M,N))$$
Here $e$ is any sufficiently large integer. The theta function was
first introduced by Hochster in his study of the direct
summand conjecture (\cite{Ho1}). We recall the basic properties of
$\theta(M,N)$ and prove a key technical result in Proposition \ref{rg1}: {the vanishing of
$\theta^R(M,N)$ implies rigidity of (M,N)}.

In section \ref{3} we study the vanishing of $\theta^R(M,N)$ when
$R$ is a hypersurface with isolated singularity. The key advantage
with such rings is that $\theta^R(M,N)$ is always defined for a
pair $M,N$, so we can ``move" the modules within the Grothendieck group into more favorable
positions  where vanishing of $\theta$ is more evident. We give a
fairly complete picture when the dimension of the ring is at most
$4$  (see \ref{dim1} to
\ref{dim2,3}). Also, Theorem \ref{moving} states that when $R$
contains a field and $\dim M + \dim N \leq \dim R$, then
$\theta^R(M,N)=0$. Our results suggest a
Conjecture (\ref{DaoConj}) that $\theta^R(M,N)$ should always vanish if $\dim R$
is {\it even}.

Seection \ref{4} focuses on rigidity over hypersurfaces in
general. We prove a new criterion for rigidity (Theorem
\ref{rig1.1}), as well as a connection to decency when one of the
modules is Cohen-Macaulay (Theorem \ref{rigidandproper}).

In section \ref{5} we  apply the results from previous sections to
give alternative proofs or extensions of relevant results in the literature. 
For example, we investigate the depth of tensor products, following the
same line of investigation done by Auslander, Huneke and
Wiegand (see \ref{vanishingiso}, \ref{HWmain}). Finally, we  give numerous examples to illustrate our
results throughout the paper and show that they are optimal in
certain senses.


After this preprint appeared, a few works 
have appeared focusing on extensions and applications of the ideas and results in here. 
For example, \cite{Da2} uses the $\Tor$-rigidity results here to extend Auslander's Theorems on 
$\Hom(M,M)$ over regular local rings to hypersurfaces, and \cite{Da3} applies such results on understanding 
Van den Bergh's definition of non-commutative crepant resolutions. Papers \cite{CeDa, Ce, Da} deal with various generalizations  to complete intersections  as well as analogous results for vanishing of $\Ext$ modules. A new interesting development is \cite{MPSW} where the Conjecture \ref{DaoConj} is studied and settled in the graded, characteristic $0$ situation using sophisticated geometric techniques.

Part of this paper was included in the author's PhD thesis at the
University of Michigan. The author would like to thank his
advisor, Melvin Hochster, for numerous invaluable discussions and
encouragements. It is a pleasure to thank Luchezar
Avramov, Ragnar-Olaf Buchweitz, William Fulton, Craig Huneke,
Mircea Musta\c t\v a, Paul Roberts and Roger Wiegand for their
patience with the author's questions regarding this work.
Special thanks must go to Greg Piepmeyer for his careful reading
of an earlier version and an anonymous referee for a very thorough 
report which vastly improves the paper. 
\section{Notation and preliminary results}\label{2}

Unless otherwise specified, all rings are Noetherian, commutative and local, and all modules are finitely generated.
A ring $(R,m,k)$ is a \textit{hypersurface} if its completion $\hat R$ has the form
$T/(f)$,
where $T$ is a regular local ring and  $f$ is in the  maximal ideal of $T$. We say
that $T$ is \textit{admissible} (as a regular local ring) if it is a power series ring over a field or a discrete
valuation ring.
If $T$ is admissible we also say that $R$ is admissible (as a hypersurface). Note that
an admissible hypersurface may be a ramified regular local ring, and thus not admissible as a regular local ring.

For a ring $R$ and a non-negative integer $i$, we set $X^iR := \{p\in \Spec(R)| \dim(R_p)\leq i\}$.
We denote by $Y(R)$ the set $X^{\dim(R)-1}$, the punctured spectrum of $R$. We denote by $G(R)$ the
Grothendieck group of finitely generated modules over $R$ and by $\overline{G}(R):= G(R)/[R]$, the reduced Grothendieck group.
Also, we let $\Sing(R) := \{p\in \Spec(R)| R_p \ \text{is not regular} \}$ be the singular locus
of $R$. For an abelian group $G$, we let $G_{\mathbb{Q}} = G\tensor_{\mathbb{Z}}\mathbb{Q}$.

Let the torsion submodule of $M$, $t(M)$, be the kernel of the map
$M\to K\tensor_RM$, where $K$ is the total quotient ring of $R$. The module $M$ is \textit{torsion} provided $t(M)=M$
and \textit{torsion-free} provided $t(M)=0$. Let $M^{*} := \Hom(M,R)$ be the dual of $M$. The module $M$ is called
\textit{reflexive} provided the natural map $M\to M^{**}$ is an isomorphism.
The module $M$ is called \textit{maximal Cohen-Macaulay}
if $\depth_RM = \dim R$. The module $M$ is said to have \textit{constant rank} if there is an integer $n$ such that  $M_p$ is $R_p$-free of rank $n$ for all minimal primes $p$ of $R$.


A pair of modules $(M,N)$ is called {\textit{rigid}} if for any integer $i\geq0$,
$\Tor_i^R(M,N)=0$ implies $\Tor_j^R(M,N)=0$ for all $j\geq i$.

One defines the finite length index of the pair $(M,N)$ as :
$$ f_R(M,N) := \min\{ i |\ l(\Tor_j^T(M,N))<\infty \ \text{for $j \geq i$} \} $$

If $M,N$ are modules over a regular local ring $T$, then for any integer $i\geq0$ such that
$f_T(M,N) \leq i$, we can define :
$$\chi_i^T(M,N) = \sum_{j\geq i} (-1)^{j-i} l(\Tor_j^T(M,N)) $$
When $i=0$ we simply write $\chi^T(M,N)$ or $\chi(M,N)$.
Serre (\cite{Se}) introduced $\chi^R(M,N)$ forty years ago as a homological definition of
intersection multiplicity for modules over a regular local ring and showed that
it satisfied many of the properties which should hold for intersection multiplicities:
\begin{thm}(Serre)\label{Serre}
Let $T$ be a regular local ring such that $\hat T$ is admissible.
Then for any pair of $T$-modules $M,N$ such that
$l(M\tensor_TN)<\infty$, we have:
\begin{enumerate}
\item
$\dim(M)+\dim(N) \leq \dim(T)$ (in other words, all modules are decent).
\item (Vanishing)
If $\dim(M)+\dim(N)<\dim(T)$, then $\chi^{T}(M,N)=0$.
\item (Nonnegativity)
It is always true that $\chi^{T}(M,N)\ge 0$.
\item (Positivity)
If $\dim(M)+\dim(N)=\dim(T)$, then $\chi^{T}(M,N)>0$.
\end{enumerate}
\end{thm}

The following ``long exact sequence for change of rings" plays a vital role in our proof of rigidity criterion
(\ref{rg1}). It follows from the famous Cartan-Eilenberg spectral sequence (\cite{Av2}, 3.3.2) and has been used frequently in 
previous work (see for example \cite{HW1}).

\begin{prop}\label{longexact}
Let $R=T/f$ where $f$ is a nonzerodivisor on $T$, and let
$M,N$ be $R$-modules. Then we have the long exact sequence of
$\Tor$s :

$$        \begin{array}{ll}
...\to \Tor_{n}^R(M,N) \to \Tor_{n+1}^T(M,N) \to \Tor_{n+1}^R(M,N)\\
\to  \Tor_{n-1}^R(M,N) \to \Tor_{n}^T(M,N) \to \Tor_{n}^R(M,N)\\
\to ...                                                         \\
\to \Tor_{0}^R(M,N) \to \Tor_{1}^T(M,N) \to \Tor_{1}^R(M,N) \to 0

\end{array}
$$\\
\end{prop}

\textbf{The infinite projective dimension locus.}
\begin{dfn}
Let $M$ be an $R$-module. One can  define the infinite projective
dimension locus of $M$ as $IPD(M):= \{ p\in \Spec(R) |\
\pd_{R_p}M_p =\infty \}$.
\end{dfn}
We gather some properties of this locus:

\begin{lem}
Let $R$ be a local hypersurface of dimension $d$. Let $M,N$ be $R$-modules.
Let $\Supp_e(M,N)=\Supp(\Tor_{2d+2}^R(M,N))$ and
$\Supp_o(M,N)=\Supp(\Tor_{2d+1}^R(M,N))$. Then we have :
\begin{enumerate}
\item  $IPD(M)$ is a Zariski closed subset of $\Spec(R)$.
\item $IPD(M)\subseteq \Supp(M)\cap \Sing(R)$.
\item $\Supp_e(M,N) \cup \Supp_o(M,N) = IPD(M)\cap IPD(N)$.
\end{enumerate}
\end{lem}

\begin{proof}
(1) Let $L = \syz_{d+1}^R(M)$. Let $F(M) = \{x\in R \ | \ L_{x} \ \text{is a free}
\ R_{x}\text{-module}\}$. For any prime $p$, $\pd_{R_p}M_p = \infty$ if and only if
$L_p$ is not free if and only if $p \supset F(M)$. So $IPD(M) = V(F(M))$.\\
(2) This is obvious. \\
(3) Let $p\in \Spec(R)$ and localize at $p$. Then $R_p$ is a hypersurface.
The assertion follows from the fact that  $\Tor_{2d+2}^{R_p}(M_p,N_p)
= \Tor_{2d+1}^{R_p}(M_p,N_p)=0$ if and only if one of $M_p,N_p$ has
finite projective dimension (see \cite{HW2}, Theorem 1.9).
\end{proof}

\textbf{The function $\theta^R(M,N)$.} \label{}

Let $R$ be a local hypersurface and assume $\hat R = T/(f)$ where $T$ is a regular local ring. The function $\theta^R(M,N)$ was introduced by Hochster (\cite{Ho1})
for any pair of finitely generated modules $M,N$ such that $f_R(M,N)<\infty$  as:
$$ \theta^R(M,N) = l(\Tor_{2e+2}^R(M,N)) - l(\Tor_{2e+1}^R(M,N)) $$
where $e$ is any integer $\geq (d+2)/2$. It is well known (see \cite{Ei}) that $\Tor^R(M,N)$ is periodic of
period 2 after $d+1$ spots, so this function is well-defined. The theta function satisfies
the following properties. First, if $M\tensor_RN$ has finite length, then:
$$\theta^R(M,N) = \chi^T(M,N)$$
As a consequence of this fact and \ref{Serre}, we have the following result from \cite{Ho1}:
\begin{prop}\label{decent}(Hochster)
Let $R$ be an admissible hypersurface and $M,N$ be $R$-modules such that $M\tensor_RN$ has finite length.
Then $(M,N)$ intersect decently if and only if $\theta^R(M,N)=0$.
\end{prop}
Secondly, $\theta^R(M,N)$ is biadditive on short exact sequences,
assuming it is defined. Specifically, for any short exact
sequence:
$$0 \to N_1 \to N_2 \to N_3 \to 0$$
and any module $M$ such that $f_R(M,N_i)<\infty$ for all
$i=1,2,3$, we have $\theta^R(M,N_2) = \theta^R(M,N_1) +
\theta^R(M,N_3)$. Similarly, $\theta(M,N)$ is additive on the
first variable. Hochster exploited these properties to give a
sufficient condition for a cyclic module in $R$ to  intersect
decently. We record Hochster's result here, in a slightly more general
form (his result was stated in terms of ideals, but the proof is virtually unchanged):

\begin{thm}(\cite{Ho1}, Theorem 0.1)\label{hochster}
Let $R$ be an admissible local hypersurface and $M,N$ be $R$-modules. Assume that :
\begin{enumerate}
\item $\Supp(M) \supseteq \Sing(R)$.
\item $[M] = 0$ in $\overline G(R)_{\mathbb{Q}}$.
\item $l(M \tensor_R N) <\infty$.
\end{enumerate}
Then $\theta^R(M,N)=0$ and $\dim M + \dim N \leq \dim R$ (In other words, any module $M$ satisfying (1) and (2) is decent).
\end{thm}

\textbf{The depth formula.}

A result by Huneke and Wiegand showed that when  all the positive $\Tor$ modules vanish,
the depths of the modules satisfy a remarkable equation:

\begin{prop}(\cite{HW1}, 2.5)\label{lem1} Let $R$ be a complete intersection. Let $M,N$ be non-zero
finitely generated modules over $R$ such that $\Tor_i^R(M,N)=0$ for all $i\geq 1$. Then:
$$ \depth(M) + \depth(N) = \depth(R)+\depth(M\tensor_RN) $$

\end{prop}

\textbf{Vanishing of $\theta^R(M,N)$ implies rigidity.}

The main technical result of this section says that, when $\theta^R(M,N)$ can be
defined and vanishes, then $(M,N)$ is rigid:

\begin{prop}\label{rg1}
Let $R$ be an admissible hypersurface and $M,N$ be $R$-modules such that $f_R(M,N)<\infty$ (so that $\theta^R(M,N)$
can be defined). Assume $\theta^R(M,N)=0$. Then $(M,N)$ is rigid.
\end{prop}

Our main tool is a celebrated Theorem  first proved by Lichtenbaum (\cite{Li}) except in a few cases. Those cases
were completed by Hochster (\cite{Ho2}):

\begin{thm}(Lichtenbaum, Hochster)\label{parchi}
Consider finitely generated modules $M,N$ over an admissible regular local ring $T$ and
an integer $i$ such that $f_T(M,N)\leq i$ (so that $l(\Tor_j^T(M,N))<\infty$ for $j \geq i$).
Then
$\chi_i^T(M,N)\geq 0$ and it is $0$ if and only if $\Tor_j^T(M,N) = 0$ for all $j\geq i$.
\end{thm}

In order to prove our rigidity result we will first need to prove
a pivotal case, when all the relevant $\Tor$s have finite length,
so we can apply Theorem \ref{parchi}.

\begin{lem}\label{pivotal}
Let $R,M,N$ be as in \ref{rg1}. Let $i$ be an integer such that $i
\geq f_R(M,N)$. Assume that $\theta^R(M,N)=0$ and $\Tor_i^R(M,N) =
0$. Then $\Tor_j^R(M,N) = 0$ for all $j\geq i$.
\end{lem}

\begin{proof}
By completion we may assume $R=T/(f)$ where $T$ is an admissible
regular local ring. We truncate the change of rings long exact
sequence for $\Tor$s (\ref{longexact}) as follows (note that
all $\Tor^T(M,N)$ vanish after $d+1$ spots):\\
$$        \begin{array}{ll}

                                   0 \to \Tor^R_{2e+2}(M,N) \\
\to \Tor_{2e}^R(M,N) \to \Tor_{2e+1}^T(M,N) \to \Tor_{2e+1}^R(M,N)\\
\to  ...                                                          \\
\to \Tor_{i}^R(M,N) \to \Tor_{i+1}^T(M,N) \to \Tor_{i+1}^R(M,N) \to C \to 0
\end{array}
$$\\
It is easy to see that all the modules in this sequence have finite lengths. Therefore we can
take the alternating sum of the lengths and get :
$$ l(C) + \chi_{i+1}^T(M,N) = (-1)^{2e+2-i}\theta^R(M,N) + l(\Tor_i^R(M,N))  = 0$$
This equation and Theorem \ref{parchi} forces $C=0$ and $ \Tor_j^T(M,N) =0$ for all $j\geq i+1$.
The conclusion of the Lemma now follows easily.

\end{proof}

Now we can prove our rigidity result :

\begin{proof} (of \ref{rg1})
We use induction on $d = \dim R$. If $d=0$, $M,N$ both have finite
length, so the previous Lemma applies. Now assume $d \geq 1$ and
$\Tor_i^R(M,N)=0$. Localizing at any $p\in Y(R)$, the punctured
spectrum of $R$, and using the induction hypothesis (note that
$R_p$ is at worst a hypersurface with dimension less than d, and
$\Tor_j^{R_p}(M_p,N_p)=0$ for $j\geq f_R(M,N)$), we may conclude
that $f_R(M,N) \leq i$. Again Lemma \ref{pivotal} can be applied
to finish the proof.
\end{proof}

\section{Hypersurfaces with isolated singularity}\label{3}

In this section we investigate the vanishing of $\theta^R$
when $R$ is a local hypersurface with isolated singularity. In this situation
$\theta^R(M,N)$ is defined for all pairs $(M,N)$ (since all higher
$\Tor$ modules have finite length). In this situation $\theta^R$
defines a bilinear map from $G(R)_{\mathbb{Q}}\times
G(R)_{\mathbb{Q}}$ to $\mathbb{Q}$, hence by Theorem \ref{rg1} it
vanishes whenever one of the modules is equivalent to $0$ in $\overline
G(R)_{\mathbb{Q}} = G(R)_{\mathbb{Q}}/{\mathbb{Q}}[R]$ (since
$\theta^R(R,-) = 0$). We record this here for convenience.
\begin{cor}
Let $R$ be an admissible hypersurface with isolated singularity.
Then $\theta^R(M,N)$ is always defined and vanishes if $M=0$ in
$\overline G(R)_{\mathbb{Q}}$.
\end{cor}

However, our investigation shall show that there are many more
cases when $\theta^R$ vanishes. Our methods and inspirations come
mostly from intersection theory. One key point is that we can
often ``move" the modules in to favorable position to show
vanishing of $\theta^R$. Since moving in the Grothendieck group is
much harder than in the Chow group, we need to make use of the
Riemann-Roch map between the two groups.

We first review some facts about Chow groups, taken from \cite{Fu} and \cite{Ro}. Let $X$ be a
Noetherian scheme.
Let $Z_iX$ be the free abelian group on the $i$-dimensional subvarieties
(integral, closed subschemes) of $X$. For any $i+1$-dimensional
subscheme $W$ of $X$, and a rational function $f$ on $W$, we can
define an element of $Z_iX$ as follows:
$$ [\Div(f,W)] = \sum_{V} \ord_V(f)[V] $$
summing over all codimension one subvarieties $V$ of $W$. Then the $i$-Chow group $\CH_i(X)$ is defined
as the quotient of $Z_iX$ by the subgroup generated by all elements of the form $[\Div(f,W)]$.
Let $\CH_*(X) = \oplus \CH_i(X)$ and $\CH_*(X)_{\mathbb{Q}} = \CH_*(X)\tensor_{\mathbb{Z}}{\mathbb{Q}}$.
An $i$-cycle (resp. cycle class) is an element in $Z_i(X)$ (resp. $\CH_i(X)$)
(by a slight abuse of notation, we also talk about a cycle as an element
of $Z_i(X)_{\mathbb{Q}}$). When $X= \Spec(R)$, where $R$ is a ring, we simply write
$\CH_*(R)$. We write $\CH^i(X)$ for $\CH_{d-i}(X)$, here $d =\dim X$.

If $R$ is local and is a homomorphic image of a regular local ring $T$ we have the important notion of
{\it{Todd class}}.  For any $R$-module $M$, let $F_*$ be the minimal free resolution of $M$ over $T$.
The {\it{Todd class}} of $M$ is defined as:
$$ \tau_{R/T}(M) := \ch(F_*)([R]) $$
Here $\ch()$ denotes the local Chern character. For much more
details about the definition and properties of the Todd class, we
refer to [Fu] or [Ro]. The Todd class has been very useful to
prove such results as  Serre's Vanishing Conjecture and the New
Intersection Theorem ([Ro]). The Todd class gives an isomorphism
of $\mathbb{Q}$-vector spaces:
$$ \tau_{R/T} : G(R)_{\mathbb{Q}}\to \CH_*(R)_{\mathbb{Q}}$$
When $R,T$ are clear we will simply write $\tau$ for $\tau_{R/T}$.
Recall that the Todd class satisfies the top terms property:
$$ \tau(M) = \sum_{\dim R/p =\dim M} l(M_p)[R/p] + \text{terms of lower dimension} $$

We collect below a number of facts that will be used frequently:

\begin{prop}\label{factschow} Let $R$ be a local ring. 
Let $M$ be an $R$-module and $d=\dim R$.
\begin{enumerate}
\item If $d>0$ and $l(M)<\infty$ then $[M]=0$ in $G(R)_\mathbb{Q}$.
\item If $R$ is regular, $\CH^i(R)_\mathbb{Q} = 0$ for $0 < i \leq d$ and $\CH^0(R)_\mathbb{Q} = \mathbb{Q}$.
\item Assume that $R$ is a homomorphic image of a regular local ring.For any $P\in \Spec R$ we have $\tau^{-1}([R/P]) = [R/P] +  \text{terms of lower dimension} $.
\end{enumerate}
\end{prop}

\begin{proof}
\begin{enumerate}
\item It is easy to see that  $[M] = l(M)[R/m]$ in $G(R)$. But take any $P\in \Spec R$
such that $\dim R/P =1$ and take $x \in m-P$, from the short exact sequence 
$$ 0 \to R/P \to R/P \to R/(P+(x)) \to 0$$
we then have $[R/(P+(x))] =0$ in $G(R)$. It follows that 
$[R/m]=0$ in $G(R)_\mathbb{Q}$, which gives what we want. 
\item Since $R$ is regular, $G(R) = \mathbb{Z}[R]$ since any module has a finite resolution by
free modules. The claim follows because $\tau([R]) = [R]$.
\item This follows immediately from the top term property.
\end{enumerate}
\end{proof}

Throughout the rest of this section we will assume that $R$ is a
local hypersurface with isolated singularity. 

\begin{prop}\label{dim1}
Assume $\dim R = 1$. Then $\theta^R(M,N)$ = 0 for all $N$ if and
only if $M$ has constant rank,or equivalently, $[M]=0$ in
$\overline{G}(R)$. 
\end{prop}

\begin{proof}
Let $p_1,p_2,...,p_n$ be the minimal primes of $R$.
Then $G(R)_{\mathbb{Q}}$ has a basis consisting of 
the classes $[R/p_1],...,[R/p_n]$ by \ref{factschow}. In
particular, since $R$ has dimension 1 and is reduced, $[R] =
\sum_{1}^{n}[R/p_i]$. Let $\alpha_{ij} = \theta^R(R/p_i,R/p_j)$.
For $i \neq j$, $p_i + p_j$ is $m$-primary, and it is not hard
(using the resolution of $R/p_i$ or \ref{ex2}) to see that $\alpha_{ij} = l(R/(p_i+p_j)) >0$.
Since $\theta^R(R,R/p_i) = 0$, we must have $\alpha_{ii} =
-\sum_{j\neq i} \alpha_{ij}$. Now, for a  module $M$, let $[M] =
\sum a_i[R/p_i]$, here $a_i$ is the rank of $M_{p_i}$. If  $a_1 =
a_2 = ... =a_n$, then $[M]=a_1[R]$, so $\theta^R(M,N) =
a_1\theta^R(R,N)= 0$. For the other direction, without loss of
generality we may assume that $a_1$ is the largest among the
$a_i$'s. Then since $0 = \theta^R(M,R/p_1) = \sum_{i=2}^{n}
\alpha_{1i}(a_i-a_1)$, we must have $a_i=a_1$ for all $i$.
\end{proof}

\begin{thm}
Assume $\dim R =d > 1$ and  $M,N$ are $R$-modules. Then $\theta^R(M,N)=0$ if
$ \dim M \leq 1$.

\end{thm}

\begin{proof}
Without affecting relevant issues we may complete and assume $R$ is a 
homomorphic image of a regular local ring. Since any module has a filtration by prime cyclic modules, we may
assume that $M = R/P$ and $N=R/Q$ for some $P,Q \in \Spec R$. If
$\dim R/P = 0$, so $P=m$, then $[R/P] = 0$ in $G(R)_\mathbb{Q}$,
and $\theta$ vanishes. Also, we may assume $Q \neq 0$. If $Q$ is
not contained in $P$, then $l(R/(P+Q)) <\infty$ because $\dim R/P
=1$, and since $\dim R/P +\dim R/Q \leq \dim R$ we have
$\theta(R/P,R/Q) = 0$ by \ref{decent}. So now we only need to
consider the case $0 \neq Q \subset P$. We claim that there is
cycle $\alpha = \sum l_i[R/Q_i] \in \CH^*(R)_\mathbb{Q}$ such that
$\alpha = [R/Q]$ and $Q_i \nsubseteq P$. Consider the element
$[R_P/Q] \in \CH^*(R_P)_\mathbb{Q}$. Since $R_P$ is regular,
$[R_P/Q] = 0$. Therefore, formally, we have a collection of primes
$q_i$ and elements $f_i$ and integers $n, n_i$ such that $n[R_P/Q]
= \sum \Div(R_P/q_i,f_i)$. Now in $\CH^*(R)_\mathbb{Q}$ we will
have $\sum \Div(R/q_i,f_i) = n[R/Q] + \sum n_i[R/Q_i]$, with $Q_i
\nsubseteq P$, which proves our claim. The fact that $[R/Q] =
\sum_i l_i[R/Q_i]$ in $\CH^*(R)_\mathbb{Q}$ means that in
$G(R)_\mathbb{Q}$, $[R/Q]=\sum_i l_i[R/Q_i]+ \text{terms of lower
dimension}$ (by part (3) of the Proposition \ref{factschow}). By
the argument at the beginning of the proof,
$\theta^R(R/P,R/Q_i)=0$ and we may conclude our proof by induction
on $\dim R/Q$.
\end{proof}

\begin{rmk}
The above argument is true whenever $R_P$ is regular. This argument could be thought of as an algebraic
``moving Lemma" for $\Spec(R)$.
\end{rmk}

When the hypersurface $R$ contains a field, we can of course apply the real
moving Lemma.

\begin{thm}\label{moving}
Suppose that $(R,m,k)$ is an excellent local hypersurface containing a field with isolated singularity. Then
$\theta^R(R/P,R/Q)=0$ if $\dim R/P + \dim R/Q \leq \dim R$.
\end{thm}

\begin{proof}
We first make some reductions. Assume we have a
counterexample on $R$. We can first make a faithfully flat
extension to replace $k$ by an algebraically closed field and then
complete
to get to the case of $R = k[[x_0,...x_d]]/(f)$, and $k$ is
algebraically closed. The condition of isolated singularity is
preserved by faithfully flat extension (cf. Lemma 2.7, \cite{Wi}).
Then by a Theorem of Artin (see \cite{CS}, Theorem 1.6), $R =
\hat{S}$, where $S$ is local hypersurface with isolated
singularity, and is essentially of finite type over $k$. But we
can descend our example to $S^h$, the Henselization of $S$, by 
standard arguments (see \cite{Ho3} or \cite{Du1}). The only issue
is how to descend the resolutions of the modules, which may be
infinite. However, note that since the resolutions of our modules
are eventually periodic, we only need to keep track of a finite
part. Then since $S^h = \underrightarrow{\lim} S_n$, where each
$S_n$ is an \'etale neighborhood of $S$, we have a counterexample
in some $S_n$. Thus we may assume $R$ is essentially of finite
type over an algebraically closed field $k$. Let's say $R = A_m$,
where $A$ is a finite $k$-algebra and $m$ is a maximal ideal in $A$. Since
$R$ has isolated singularity, we have $\Sing(A)$ is a disjoint
union of $\{m\}$ and some closed subset $Y \subset \Spec A$.

Let $X = \Spec(A) -\{m\}-Y$. Then $X$ is a smooth quasi-projective
variety, so by the Moving Lemma (see \cite{Fu},11.4) one can find
a cycle $\alpha = \sum_i n_i V(Q_i)$ in $\CH^*(X)$ such that for
each $i$, $\dim (V(Q_i)\cup V(P))\leq 0$, that is the intersection
consists only of points in $X$. When we restrict all the cycles
and divisors to $\Spec(R)$ we will have $[R/Q] = \sum_j [R/Q_j]$
in $\CH^*(R)$ and for all $j$, $V(P)\cup V(Q_i) \subset \{m\}$. In
$G(R)_{\mathbb{Q}}$ this means $[R/Q] = \sum_j [R/Q_j] +
\text{terms of lower dimension}$. Since $l(R/P+Q_j) <\infty$ and
$\dim R/P + \dim R/Q_j \leq \dim R$ we have
$\theta^R(R/P,R/Q_j)=0$ by \ref{decent} and an induction on $\dim
R/Q$ finishes the proof.
\end{proof}

\begin{cor}\label{dim2,3}
If $\dim R=2$ or  $\dim R=3 $ and $\CH^1(R)_{\mathbb{Q}} = 0$ or $\dim R =4$ and $R$ excellent and contains a field, then
$\theta^R(M,N)=0$ for all pairs $(M,N)$.
\end{cor}

\begin{proof}
It suffices to assume that $M,N$ are cyclic prime modules, let's
say $M=R/P,N=R/Q$. Then by the previous Theorems we only need to
worry if both of them have dimension at least $2$. If $\dim R=2$, they must both be  $R$ (note that since $R$ is normal and local, it is a domain), 
thus $\theta^R$ certainly vanishes. If $\dim R=3$, assume that $\dim R/P=2$ (otherwise $R/P$ would be $R$).
We can then complete $R$ without loss of generality. Then part (3) of \ref{factschow} shows that in $G(R)$, $[R/P]$ is
equal to a formal sum of cyclic primes of dimension $\leq 1$,
forcing $\theta^R(R/P,R/Q)$ to be $0$. Finally, if $\dim R=4$ and
$R$ contains a field we can apply \ref{moving} and assume $\dim
R/P + \dim R/Q \geq 5$. Then one of the primes, say $P$ is height
$1$ (if the minimal height is $0$ the assertion is trivial). We will be done if we can show
that $\CH^1(R)=0$. But by the Grothendieck-Lefschetz Theorem, the
Picard group of $X = \Spec(R)-\{m\}$ is $0$. Since $X$ is
regular, the Picard group of $X$ is the same as $\CH^1(X)=
\CH^1(R)$.
\end{proof}

\begin{eg}
Let $R$ be a dimension $3$ hypersurface which is an UFD, for example $R=\mathbb C[[x,y,u,v]]/(xy+f(u,v))$ where $f\in \mathbb C[[u,v]]$ is irreducible.  Then $\CH^1(R)=0$ and $\theta$ vanishes on any pair of
modules over $R$.
\end{eg}

Our next result is an algebraic Bertini type Theorem for the
vanishing of $\theta^R$. It is most useful when $M=N$ (so in a
sense when moving them apart is the hardest).

\begin{thm}
Assume that $\dim R \geq 2$. Suppose $M,N$ are $R$-modules such that there is an element
$x \in \Ann{M} \cap \Ann(N)$ such that $R/(x)$ is still a
hypersurface with isolated singularity. Then $\theta^R(M,N)=0$.
\end{thm}

\begin{proof}
Let $L = \syz_R^{2d}(N)$. Then $\theta^R (M,N) = \theta^R(M,L)$. Also, $L$ is maximal Cohen-Macaulay: in
particular, $x$ is a nonzerodivisor on $L$. Thus $\Tor_i^R(M,L) = \Tor_i^{R/x}(M,L/(x))$, so it is enough
to prove that $\theta^{R/(x)}(M,L/(x)) = 0$ (the assumption that $R/(x)$ is still a hypersurface
with isolated singularity ensures that $\theta^{R/(x)}$ is well-defined). We define a map
$\alpha : G(R) \to G(R/(x))$ as follows:
$$ \alpha([M]) = [\Tor_0^R(M,R/(x))] - [\Tor_1^R(M,R/(x))] $$
We need to show that $\alpha$ is well-defined. The only thing that
needs to be checked is that if $0 \to M_1 \to M_2 \to M_3 \to 0$
is a short exact sequence of $R$-modules, then $\alpha([M_2]) =
\alpha([M_1]) + \alpha([M_3])$. But this follows from tensoring
the exact sequence with $R/(x)$, and note that $\Tor_i^R(M,R/(x))
=0$ for all $R$-modules $M$ and $i>1$.

It is easy to see that  $\alpha([N]) = [0]$ and $\alpha([L]) = [L/(x)]$
(since $x$ kills $N$ and is $L$-regular, respectively). But since $L = \syz^{2d}(N)$ we have
$[L] = [N] + n[R]$ for some integer $n$. Applying $\alpha$ we get:
$ [L/(x)] = n[R/(x)]$ in $G(R/(x))$. But then $\theta^{R/(x)}(M,L/xL) = 0$, which is what we want.
\end{proof}

Now we will consider a special case, when $R$ is a local ring of
an affine cone of a  projective variety. That is, $R=A_{\mathfrak
m}$, where $A$ is a graded hypersurface over a field $k$ whose
homogeneous maximal ideal is $\mathfrak m$. The following result
by Kurano is very helpful (cf. \cite{Ku1}, Theorem 1.3).
\begin{thm}(Kurano)\label{kurano}
Let $A,R$ be as above. Let $X = \Proj(A)$. Assume that $R$ has
isolated singularity.Then $X$ is smooth and $\CH^*(X)$ becomes a
(graded) commutative ring with the intersection product. Let $h
\in \CH^1(X)$ represent the hyperplane section (alternatively, the
first Chern class of the invertible sheaf ${\mathcal O}_X(1)$).
Then there is a graded isomorphism of ${\mathbb{Q}}$-vector spaces
from $\CH^*(X)/h\CH^*(X)$ to $\CH^*(R)$.
\end{thm}

\begin{prop}
Let $X,R$ as above and $i\geq 0$ an integer. If $\CH^{i+1}(X)_{\mathbb{Q}} = \mathbb{Q}$
then $\CH^{i+1}(R)_{\mathbb{Q}} = 0$.
\end{prop}

\begin{proof}
We only need to observe that in $\CH^*(X)_{\mathbb{Q}}$, multiplication by $h$ is a nonzero map from
$\CH^i(X)_{\mathbb{Q}}$ to $\CH^{i+1}(X)_{\mathbb{Q}}$ (just look at $h.h^{i}=h^{i+1}$). Hence if
$\CH^{i+1}(X)_\mathbb{Q} = \mathbb{Q}$
then $ \CH^{i+1}(X)_\mathbb{Q} =h\CH^{i}(X)_\mathbb{Q}$, so $\CH^{i+1}(R)_\mathbb{Q}=0$ by Kurano's
result.
\end{proof}

This allow us to exploit many results in the literature about the Chow groups of projective varieties.
For example, we have the following:

\begin{thm}(\cite{ELV},2.3) Let $X \subset \mathbb{P}^d_k$ be an irreducible hypersurface of degree
$s$. If:
$$ {{i+s} \choose {i+1}} \leq d-i$$
Then $\CH^{d-1-i}(X)_{\mathbb{Q}} = \mathbb{Q}$.
\end{thm}

\begin{cor}
Let $X$ be a smooth hypersurface of degree $s$ in $\mathbb{P}^d_k$. Let $R$ be the local ring
at the homogenous maximal ideal of the affine cone of $X$. Let $n$ be the biggest integer such that:
$$ {{n+s} \choose {n+1}} \leq d-n$$
Then $\theta^R(M,N) = 0$ if $\dim M \leq n+1$.
\end{cor}

\begin{proof}
By the two previous results we have $\CH^{d-1-i}(R)_{\mathbb{Q}} = 0$ for $i\leq n$, in other words,
$\CH_{i}(R)_{\mathbb{Q}} = 0$ for $i\leq n+1$. Thus $\tau([M]) = 0$, hence $[M] = 0$ in $G(R)_{\mathbb{Q}}$.
\end{proof}

\begin{rmk}
When $s=2$, we have $n+1 = \lfloor d/2 \rfloor $. This gives a strong version of \ref{moving}.
\end{rmk}

Next we want to discuss a conjecture, attributed to Hartshorne
(see \cite{Ha1}, page 142), which could be relevant to our
interest:

\begin{conj}\label{HarConj}(R.Hartshorne) Let $X$ be a smooth projective, complete intersection variety
in $\mathbb{P}^n_k$. Then $\CH^i(X)_{\mathbb{Q}}= \mathbb{Q}$  for
$i<{\dim X}/2$.
\end{conj}

It is interesting to observe that Hartshorne's conjecture,
together with \ref{moving}, shows that if $\dim R$ is even, when
$R$ is the local ring at the origin of the affine cone of $X$,
then $\theta^R$ always vanishes.

\begin{cons}(of Hartshorne's conjecture)
Let $X$ be a smooth hypersurface in $\mathbb{P}^d_k$. Assume that
$d$ is even. Let $R$ be the local ring at the origin of the affine
cone of $X$. Then $\theta^R$ always vanishes.
\end{cons}

\begin{proof}
Let $d=2n$. We have $\dim X = 2n-1$, so by Hartshorne' conjecture with $Y=\mathbb{P}^d_k$ we have
$\CH^i(X)_{\mathbb{Q}} ={\mathbb{Q}}$
for $i \leq n-1$. Thus $\CH^i(R)_{\mathbb{Q}} = 0$ for $i\leq n-1$, in other words, $\CH_i(R)_{\mathbb{Q}} = 0$
for $i \geq n+1$ . So in the Grothendieck group $G(R)_{\mathbb{Q}}$,
any module can be represented  as a sum of cyclic prime modules of dimension $\leq n$. But since
$\dim R = 2n$, for any such pair of modules $(R/P,R/Q)$  we must have $\theta^R(R/P,R/Q)=0$ by \ref{moving}.
\end{proof}

In view of this and our knowledge of dimensions $2$ and $4$, we feel it is reasonable to make:

\begin{conj}\label{DaoConj}
Let $R$ be a hypersurface with isolated singularity. Assume that $\dim R$ is even and $R$ contains a field.
Then $\theta^R$ always vanishes.
\end{conj}

We observe that the values of $\theta^R$ only depend on its values on pairs of
maximal Cohen-Macaulay (MCM) modules (as one can replace the
modules by their high syzygies). It is worth noting that the theta
function is closely related to the notion of ``Herbrand
difference" defined on a pair of MCM modules using stable cohomology by Buchweitz in
\cite{Bu}. In fact, for MCM modules $M,N$ we have $\theta^R(M,N)$ and the Herbrand difference $\text{h}(M,N^*)$ agree up to sign (we thank Ragnar-Olaf Buchweitz for explaining this connection to us).

Recall that  a complete local hypersurface $R$ is called a \textit{simple singularity} if it is isomorphic to $T/(f)$, where $T=k[[x_0,x_1,...,x_d]]$ for some  $d>0$,
$k$ an algebraically closed field of characteristic $0$, and $f$ has  one of the following forms :\\
$ (A_n)\ \ \ \ x_0^2+x_1^{n+1} + x_2^2+...+x_d^2 \ \ \ \ \ \ (n\geq 1)$\\
$ (D_n)\ \ \ \ x_0^2x_1 + x_1^{n-1} +x_2^2+...+x_d^2 \ \ \ \ (n\geq 4)$\\
$ (E_6)\ \ \ \ x_0^3 + x_1^4 +x_2^2+...+x_d^2 \ \ \ \ $\\
$ (E_7)\ \ \ \ x_0^3 + x_0x_1^3 +x_2^2+...+x_d^2 \ \ \ \ $\\
$ (E_8)\ \ \ \ x_0^3 + x_1^5 +x_2^2+...+x_d^2 \ \ \ \ $\\

For hypersurfaces, simple singularity is the same as finite representation type, that is, the
group of isomorphism classes of indecomposable MCM modules is finite. The Grothendieck
group of MCM modules over simple singularities have been computed completely (see \cite{Yo}, 13.10). One striking feature is that in even dimensions, all the Grothendieck groups are torsion after we kill the class of $[R]$. Thus
we have the following result, which confirm Conjecture \ref{DaoConj} (we thank the referee for pointing out that it also follow from 10.3.8 in \cite{Bu}):

\begin{cor}\label{MainConj}
Let $R$ be a hypersurface with isolated, simple singularity of even dimension. Then $\theta^R$ always vanishes.
\end{cor}

\section{Rigidity over hypersurfaces}\label{4}

By  virtue of Proposition \ref{rg1} and the results in the previous section, we have a lot
of results about rigidity of modules when the hypersurface has an isolated singularity. 
On general hypersurfaces, however, we need to be more careful about using the function $\theta^R$.
Typically, we need some extra conditions to show that $\theta$ is defined for all the modules
in the short exact sequences involved. In the last section we will give plenty of examples to show that these
conditions are unavoidable. 

We will now state two immediate corollaries of \ref{rg1}. The first appeared implicitly in the work
of Lichtenbaum, (\cite{Li}) :
\begin{cor}\label{lichten}
Let $R$ be an admissible hypersurface and $M,N$ be finitely generated
$R$-modules. If $M$ or $N$ has finite projective dimension, then $(M,N)$ is rigid.
\end{cor}

\begin{cor}\label{rg2}
Let $R$ be an admissible hypersurface and $M,N$ be finitely
generated $R$-modules. Assume that $\pd_{R_p}M_p <\infty$ for all
$p\in Y(R)$ (the punctured spectrum of $R$) and $[N]=0$  in
$\overline{G}(R)_{\mathbb{Q}}$. Then $(M,N)$ is rigid.
\end{cor}

\begin{proof}
The first assumption ensures that $f_R(M,N)<\infty$ for all $N$, hence
$\theta^R(M,N)$ can be defined for all $N$. Then the second assumption forces
 $\theta^R(M,N) = 0$.
\end{proof}
Another immediate corollary of our result is the first ``rigidity" Theorem in a
paper of Huneke and Wiegand (\cite{HW1}).

\begin{cor}\label{HWrg}
(\cite{HW1}, 2.4). Let $R$ be an admissible hypersurface and $M,N$ be
$R$-modules. Assume:
\begin{enumerate}
\item $M\tensor_RN$ has finite length.
\item $\dim(M) + \dim(N) \leq \dim(R)$.
\end{enumerate}
Then $(M,N)$ is rigid.
\end{cor}

\begin{proof}
Suppose $R=T/(f)$ where $T$ is regular local. In this case $\theta^R(M,N) =\chi^T(M,N)$,
so by Serre's vanishing Theorem, it must be $0$.
\end{proof}

The next result introduces a class of rigid modules not necessarily having finite projective
dimension. To state it, recall the definition: $IPD(M):= \{ p\in \Spec(R) |\ \pd_{R_p}M_p =\infty \}$ (section 1)
\\
\begin{thm}\label{rig1.1}
Let $R$ be an admissible hypersurface, and $M$ be an $R$-module such that $[M] =0$
in $\overline G(R)_{\mathbb{Q}}$. Assume that $IPD(M)$ is either $\emptyset$ or is equal to $\Sing(R)$. Then
$M$ is rigid.
\end{thm}

\begin{proof}

If $IPD(M)=\emptyset$ then $\pd_RM<\infty$, so there is nothing to prove.
 Assume that $IPD(M)= \Sing(R) \neq 0$.
We again use induction on $d = \dim R$. If $d=0$ the condition that $[M]=[0]$
in $\overline G(R)_{\mathbb{Q}}$ implies that $\theta^R(M,N)$ is defined and equal to $0$ for any $R$-module $N$.
Suppose $d>0$ and $\Tor_i^R(M,N)=0$ for some $i$. We localize at any prime $p\in Y(R)$.
Both conditions on $M$ localize, so by the induction hypothesis, $\Tor_j^{R_p}(M_p,N_p)=0$
for $j \geq i$. This forces either $\pd_{R_p}M_p$ or $\pd_{R_p}N_p$ to be finite (see Theorem 1.9, \cite{HW2}).
But since $I(M) = \Sing(R)$, $N$ must have finite projective dimension on $Y(R) \cap \Sing(R)$.
So $N$ has finite projective dimension on $Y(R)$, hence $\theta(M,N)=0$, finishing the proof.
\end{proof}

The following will be useful for our application to torsion of tensor products.

\begin{lem}\label{vanishinglem}
Let $R$ be an admissible hypersurface and $M,N$ be $R$-modules.
Assume that:
\begin{enumerate}
\item $\Tor_1^R(M,N)=0$
\item $\depth(N)\geq 1$ and $\depth(M\tensor_RN)\geq 1$.
\item $f_R(M,N)<\infty$.
\end{enumerate}
Then $\Tor_i^R(M,N)=0$ for $i\geq 1$.
\end{lem}

\begin{proof}
The depth assumptions ensure that we can choose $t$ a nonzerodivisor for both $N$ and $M\tensor_RN$.
Let $\overline{N} = N/tN$. Tensoring the short exact sequence :
 \[ \xymatrix {0 \ar[r] &N \ar[r]^{t} &N \ar[r] &\overline{N} \ar[r] &0} \]
with $M$ and using (1) we get :
\[ \xymatrix{0 \ar[r] &\Tor_1^R(M,\overline{N}) \ar[r] &M\tensor_RN \ar[r]^{t} &M\tensor_RN \ar[r] &M\tensor_R \overline{N} \ar[r] &0}\]
which shows that $\Tor_1^R(M,\overline{N}) = 0$. But condition (3) is satisfied for both of the pairs $(M,N)$ and $(M,\overline{N})$ and so :
$$ \theta^R(M,\overline{N}) = \theta^R(M,N) - \theta^R(M,N)= 0$$
The conclusion then follows from \ref{rg1} and Nakayama's Lemma.
\end{proof}

The next result shows some connection between rigidity, decency
and a property of modules first studied by Auslander (\cite{Au}).

\begin{thm}\label{rigidandproper}
Let $R$ be an admissible hypersurface. For a Cohen-Macaulay $R$-module $M$, the following are
equivalent:
\begin{enumerate}
\item $(M,N)$ is rigid for all $N$ such that $l(M\tensor_RN)<\infty$.
\item Every $M$-sequence is an $R$-sequence.
\item $M$ is decent ($\dim M +\dim N \leq \dim R$ for all $N$ such that $l(M\tensor_RN)<\infty$).
\item $\theta^R(M,N)=0$ for any $N$ such that $l(M\tensor_RN)<\infty$.
\end{enumerate}
\end{thm}
\begin{proof}
Assume (1). Then we can prove (2) by adapting the argument in
Auslander paper (\cite{Au},4.1) (which assumed that $(M,N)$ is
rigid for all $N$ but did not need $M$ to be Cohen-Macaulay). We
give a sketch here. Let $X$ be the free resolution of $M$. Let
$\bold x$ be a full $M$-sequence (so its length is $\dim M$ and
$l(M/(\bold {x})) <\infty$). Let $Y$ be the Koszul complex on
$\bold x$. Then the total complex $X\tensor_RY$ is acyclic.
Filtering that complex by $F_p(X\tensor Y) = \sum_{p\geq r}\sum
_qX_r\tensor_RY_q$ ,we obtain a spectral sequence with $E^2_{p,q}
= H_p(X\tensor_RH_q(Y))$. By assumption $H_n(X\tensor_RY)=0$, so
$E^{\infty}_{p,q}=0$ for $p,q>0$. We also have $E^2_{p,q}=0$ for
$p,q<0$. Hence $E^2_{1,0}=E^i_{1,0}$ for $i>1$. But
$E^{\infty}_{1,0}=0$ which implies that $H_1(X\tensor_RH_0(Y)) =
E^2_{1,0}=0$. Since $M\tensor_RH_0(Y) = M/(\bold x)$ has finite
length, we must have $0 = H_p(X\tensor_RH_0(Y)) = E^2_{p,0}$ for
all $p \geq 1$. By induction we will have:
$$ 0 = E^2_{p,q} =  H_p(X\tensor_RH_q(Y))$$
for all $p\geq 1$ and $q\geq 0$ (note that since $\bold{x}R$ kills all the modules $H_q(Y)$, $M\tensor_RH_q(Y)$
has finite length so we can apply (1)). But since $E^{\infty}_{p,q}=0$ for $p,q>0$, we have:
$$ 0 = E^2_{0,q} = H_0(X\tensor_RH_q(Y)) = M\tensor_RH_q(Y)$$
for each $q>0$. This forces $H_q(Y)=0$ for $q>0$, hence $\bold x$ is an $R$-sequence. If $\bold x$ is not
a full $M$-sequence, we can always add more elements and reach the same conclusion.\\
Assume (2). Let $N$ be an $R$ module such that $l(M\tensor_RN)<\infty$. Then we can find a full system
of parameters $\bold x$ on $M$ such that $\bold x \subset \Ann(N)$. As $M$ is Cohen-Macaulay, $\bold x$
is also a full $M$-sequence. By assumption, $\bold x$ is an $R$-sequence. Thus:
$$\dim N \leq \dim R/(\bold x)= \dim R -\dim M$$
Finally, $(3)\Rightarrow (4) \Rightarrow (1)$ is just the proof of   \ref{HWrg}.
\end{proof}

This result gives necessary conditions for rigidity that are
easier to check than rigidity itself. The conditions (2) and (3)
are quite familiar. They have played a vital role in a group of
Theorems and conjectures known as the ``homological conjectures" (see \cite{Ho4}, \cite{PS}, \cite{Ro}).

\section{Applications and examples}\label{5}
In this section we apply our results on a number of topics which involve decency or rigidity. We end by giving some examples
to complement our results. \\

\textbf{Torsion on tensor products.}\\

In this part we shall apply our rigidity results to show that
tensor products rarely have good depths:
\begin{thm}\label{vanishingiso}
Let $R$ be an admissible hypersurface with dimension $d\geq2$. Let $M,N$ be $R$-modules. Assume that:\begin{enumerate}
\item $R$ has isolated singularity.
\item $M\tensor_RN$ is torsion-free.
\item $\depth_R M\tensor_RN \geq 2$.
\end{enumerate}
Then $\Tor_i^R(M,N)=0$ for $i\geq 1$.
\end{thm}

\begin{proof}
Condition (1) makes sure that $\theta^R(M,N)$ is defined for any
pair of modules $(M,N)$. We may assume that $M,N$ are torsion-free
(by the argument in the proof of 2.4 in \cite{HW1}) we
repeat it here for the reader's convenience). Now there is an exact sequence:
$$ 0 \to M \to R^{\lambda} \to M_1 \to 0$$
Here $\lambda = \lambda(M^*)$, the number of generators of $M^*$. This exact sequence is called the
\textit{pushforward} of $M$ (see \cite{HJW}).
By tensoring the pushforward exact sequence of $M$ with $N$ , we
get:
$$0 \to \Tor_1^R(M_1,N) \to M\tensor_RN \to N^{\lambda} \to M_1\tensor_RN \to 0$$
By condition (1) $N$ is generically free, so $\Tor_1^R(M_1,N)$ is torsion, and it must be $0$
since $M\tensor_RN$ is torsion-free. Since $\depth_R M\tensor_RN \geq 2$, $\depth_RN\geq1$
(since $N$ is torsion-free and $d\geq2$), we must have $\depth_R M_1\tensor_RN \geq 1$.
Now the desired assertion follows from Lemma \ref{vanishinglem}.
\end{proof}

In the dimension $2$ case, we can do a little bit better:

\begin{prop}\label{dim2normal}
Let $R$ be an admissible hypersurface of dimension $2$. Assume further that $R$ is  normal. Let
$M,N$ be $R$ modules such that $M\tensor_R N$ is torsion-free. Then $\Tor_i^R(M,N)=0$ for $i\geq 1$ and $M$ or $N$ has finite projective dimension.
\end{prop}

\begin{proof}
We may assume $M$ is torsion-free. Now, let $M_1$ be the pushforward of $M$. We have
$\Tor_1^R(M_1,N)=0$. By the fact that $R$ is an isolated singularity of dimension 2 and \ref{dim2,3}, every
module is rigid, so $\Tor_i^R(M_1,N)=0$ for $i>1$, which gives the desired conclusion (cf. 1.9 of \cite{HW1}).
\end{proof}

\begin{rmk}
It was asked in \cite{HW1} (4.1 and the discussion before 5.3) whether or not there are two non-free
reflexive modules over a hypersurface of dimension $2$ such that their tensor product is torsion-free. In general, such pairs of modules exist. For example, let $R=k[[x,y,z]]/(xy)$ and $M=N=R/(x)$. But
with the extra assumptions of the above, such modules can not exist. For by the conclusion, one of
them must have finite projective dimension, and, being maximal Cohen-Macaulay, must be free.
\end{rmk}

To illustrate the efficiency of using $\theta^R$ for rigidity, we will give a  short proof of
one of the key results of \cite{HW1}:

\begin{thm}\label{HWmain}(\cite{HW1},2.7)
Let $R$ be an admissible hypersurface and $M,N$ be $R$-modules, at least one of which has constant rank.
If $M\tensor_RN$ is reflexive, then $\Tor_i^R(M,N)=0$ for $i>0$.
\end{thm}

\begin{proof}
We will use induction on $d = \dim R$. If $d=0$, the constant rank
condition means one of the modules must be free, and the
conclusion follows trivially.  Now assume $d\geq 1$. By the
induction hypotheses, $l(\Tor_i^R(M,N))<\infty$ for $i>0$. As in the proof of \ref{vanishingiso} we can
assume both $M,N$ are torsion free. In particular, they must have depth at least
$1$. Let $M_1$ be the pushforward of $M$:
$$ 0 \to M \to F \to M_1 \to 0$$
Then by the same reason as in proof of \ref{vanishingiso}, we have $\Tor_1^R(M_1,N)=0$. So we have:
$$ 0 \to M \tensor_RN \to F\tensor_RN \to M_1\tensor_RN \to 0$$
By the depth Lemma, we get $\depth(M_1\tensor_RN)\geq 1$. Finally,
since $l(\Tor_i^R(M,N))<\infty$ for $i>0$ we must have
$f_R(M,N)<\infty$. Applying \ref{vanishinglem} for $M_1$ and $N$,
we get $\Tor_i^R(M_1,N)=0$ for $i>1$, which implies
$\Tor_i^R(M,N)=0$ for $i>0$.
\end{proof}

\textbf{Hypersurfaces in Projective spaces.}\\

\begin{thm}\label{projhyper}
Let $k$ be a field. Let $X \subset \PP_k^n$ be a smooth
hypersurface. Let $U,V$ be subvarieties of $X$ such that $\dim U +
\dim V \geq \dim X$. Assume that $[U] = h.[U']$ in
$\CH^*(X)_{\mathbb{Q}}$, here $h$ is the hyperplane section. Then
$U\cap V \ne \emptyset$.
\end{thm}

\begin{proof}
Let $X = \Proj(A)$ where $A = k[x_0,...,x_n]/(F)$. Let $R$ be the
local ring at the origin of $A$. Suppose $P,Q \in \Spec(R)$ define
$U,V$ respectively. Our assumption becomes : $R$ is a hypersurface
with isolated singularity of dimension $n$, $\dim R/P +\dim
R/Q\geq n+1$, and $[R/P] = 0$ in $\CH^*(R)_{\mathbb{Q}}$ (by
Kurano's Theorem \ref{kurano}). We need to show that $R/P$, as a
\emph{module}, is decent. Now, by \ref{hochster} we are done if we
can show: there exist a module $M$ such that $M =0$ in $\overline
G(R)_{\mathbb{Q}}$, and $\Supp(M) = \Supp(R/P)$. We first pick $M
= R/P$. This may not guarantee that $M =0$ in $\overline
G(R)_{\mathbb{Q}}$, because by Riemann-Roch:

$$ \tau(M) = [R/P] + \sum_{i}n_i[R/p_i]  $$
Here the $p_i$s are in $\Supp(R/P)$, but have smaller dimensions.
Our strategy will be to replace them one by one by elements of
even smaller dimensions. Let's look at $p_1$. Replacing $M$ by a
multiple of $M$ if necessary, we may assume $n_1 \in \ZZ$. Next,
we replace $M$ by $ M' = {p_1}^aM\oplus (R/p_1)^{\oplus{b}}$ for
some $a,b \in \ZZ$. Then
$$ \tau(M') = \tau(M) - \tau(M/p_1^aM) + b \tau(R/p_1)$$
We now choose $a,b$ such that $b- l(M_{p_1}/{p_1}^aM_{p_1}) =
-n_1$. Then $p_1$ is replaced in the representation of $M$ by some
elements of smaller dimension. Repeating this process, we will get
to dimension $0$, which is $0$ in $\CH^*(R)_{\mathbb{Q}}$. So we
get a module $M$ such that $\tau(M) = n[R/P] = 0$, which is what
we need.
\end{proof}

\textbf{An extension of $\theta^R(M,N)$.}\\

In this subsection we extend the definition of $\theta^R(M,N)$ slightly.
For any $R$-module $L$ , we define the class of $L$ in
$Z_*(R)$ as :
$$cl(L) := \sum_{p \in \Min(L)} l_{R_p}(L_p)[R/p]  $$
and for a pair of modules $(M,N)$, let:
$$\alpha^R(M,N):= cl(\Tor_{2e+2}^R(M,N)) - cl(\Tor_{2e+1}^R(M,N))$$
Here $e$ is any integer bigger than the dimension of $R$. Note that if the $\Tor$s have finite
length then $\alpha^R(M,N)=\theta^R(M,N)[R/m_R]$.
We then have the following corollary of \ref{rg1} :
\begin{cor}
Let $R$ be an admissible hypersurface, and $M,N$ be $R$-modules. Assume that $\alpha^R(M,N)=0$ in $Z_*(R)$.
Then $(M,N)$ is rigid.
\end{cor}

\begin{proof}
We use induction on $d = \dim R$. If $d=0$ then $\alpha^R(M,N)$
coincides with $\theta^R(M,N)$. Suppose $d>0$ and
$\Tor_i^R(M,N)=0$ for some $i>0$. Then by localizing at primes on
the punctured spectrum and induction hypotheses it follows that
$l(\Tor_j^R(M,N))<\infty$ for $j \geq i$. Then the vanishing of
$\alpha$ again implies the vanishing of $\theta^R$, and we are
done by \ref{rg1}.
\end{proof}

\begin{cor}
Let $R$ be an admissible hypersurface, and $M$ be an $R$-module. Assume that the minimal
resolution of $M$ over $R$ is eventually periodic of period 1. Then $M$ is rigid.
\end{cor}
\begin{cor}
Let $R$ be an admissible hypersurface, and $M,N$ be $R$-modules. Assume that there are
integers $0<a<b$ such that $b-a$ is an odd integer and $\Tor_a^R(M,N)=\Tor_b^R(M,N)=0$.
Then $\Tor_i^R(M,N)=0$ for $i\geq a$.
\end{cor}

\begin{proof}
Let $M' = \syz_{a-1}M\oplus\syz_{b-1}M$.  Since $b-a$ is odd and the resolution of
$M$ is eventually periodic of period 2, it follows that  the resolution of
$M'$ is eventually periodic of period 1. Because
$\Tor_i^R(M',N)=\Tor_{a+i-1}^R(M,N)\oplus\Tor_{b+i-1}^R(M,N)$,
for $i\geq1$, the conclusion follows from the rigidity of $M'$.
\end{proof}

\begin{rmk}
The case $b=a+1$ is a well-known result by Murthy (\cite{Mu}), and an asymptotic version
(i.e when $a,b$ are big enough) was proved in (\cite[3.1]{Jo2}).
\end{rmk}

Finally, we will give some examples to demonstrate that many of our technical conditions
can not be removed, and some statements can not be reversed. We begin with a Lemma giving a general 
situation when rigidity fails to hold.

\begin{lem}\label{ex2}
Let $R$ be a hypersurface and $M$ a  Cohen-Macaulay $R$-module. Assume there exists a
Cohen-Macaulay $R$-module $N$ such that:
\begin{enumerate}
\item $l(M\tensor_RN)<\infty$.
\item $\dim M + \dim N = \dim R + 1$.
\end{enumerate}
Then  $\Tor_i^R(M,N)=0$ if and only if $i$ is an odd integer.
\end{lem}

\begin{proof}
Suppose $R=T/(f)$, where $T$ is regular local. Then $\dim M + \dim N = \dim T$. Since
both $M,N$ are Cohen-Macaulay, we have $\depth M + \depth N = \depth T$ as well. By 2.2 of \cite{HW2}
we have $\Tor_i^T(M,N)=0$ for all $i\geq 1$. A glance at the change of rings exact sequence
gives $\Tor_1^R(M,N)=0$ and $\Tor_{i+2}^R(M,N)\cong \Tor_i^R(M,N)$ for $i\geq0$. But then
$\Tor_2^R(M,N)\cong M\tensor_R N \neq 0$.
\end{proof}

\begin{eg}
Corollary \ref{dim2,3} shows that  $\theta^R(-,-)$ can 
vanish even when none of the modules is $0$ in
$\overline G(R)_{\mathbb{Q}}$. For example, one can take the
hypersurface $R$ to be the local ring at the ideal $(x,y,z)$ of
$A=\mathbb{C}[x,y,z]/(x^3+y^3+z^3)$, then $X=\Proj(A)$ is an
elliptic curve, therefore $\CH_*(R)_{\mathbb{Q}}$, hence
$G(R)_{\mathbb{Q}}$ and $\overline G(R)_{\mathbb{Q}}$, are infinite
dimensional $\mathbb{Q}$-vector spaces. But \ref{dim2,3} shows that $\theta$ always vanishes over $R$. 
\end{eg}

\begin{eg}
In this example we will give a module $M$ such that
$\theta^R(M,-)$ does not always vanish, but $M$ is still rigid and
decent. Let $R=k[[x,y,u,v]]/(xu-yv)$ and $P=(x,y),Q=(x,v)$. Let $M
= R/P\oplus R/P\oplus R/Q$. It is easy to check that there is an
exact sequence:
$$ 0 \to Q \to R^2 \to P \to 0$$
which shows that $R/P+R/Q=0$ in $G(R)$. So $\theta^R(M,-) =
\theta^R(R/P,-)$. Clearly $\theta^R(R/P,-)$ is not always $0$,
because $\theta^R(R/P,R/Q)=1$. It remains to show that $M$ is
rigid and decent. Let $M'=R/P\oplus R/Q$. Then by the argument
above $M'=0$ in $G(R)$, so it is rigid and decent since $R$ is an
isolated singularity. So for any module $N$, $\Tor_i^R(M,N)=0$
$\Rightarrow$ $\Tor_i^R(R/P,N)=\Tor_i^R(R/Q,N)=0$ $\Rightarrow$
$\Tor_i^R(M',N)=0$ $\Rightarrow$ $\Tor_{i+1}^R(M',N)=0$
$\Rightarrow$ $\Tor_{i+1}^R(R/P,N)=\Tor_{i+1}^R(R/Q,N)=0$
$\Rightarrow$ $\Tor_{i+1}^R(M,N)=0$. As for decency, observe that
$M$ and $M'$ have the same support, and as decency only depends on
the support, $M$ must be decent as well.
\end{eg}

\begin{eg}
Let $R=k[[x,y,u,v,t]]/(xu-yv)$ and let $M=R/(x,y,t)$. Then $M$ is not rigid (let $N=R/(u,v)$ and
use \ref{ex2}). However, the exact sequence:
 \[ \xymatrix {0 \ar[r] &R/(x,y) \ar[r]^{t} &R/(x,y) \ar[r] &\overline{M} \ar[r] &0} \]
shows that $[M]=0$ in $\overline G(R)$. It is easy to check that $\Sing(R)=V((x,y,u,v))$ and
$IPD(M) = \{(x,y,u,v,t)\} = \{m_R\}$. This example shows that the technical requirements for
rigidity in Theorem \ref{rig1.1}  can not be relaxed.
\end{eg}


\begin{eg}
Let $R=k[[x,y,u,v]]/(xu-yv)$, $M=(x,y), N=(u,y)$. Then $M\tensor_RN \cong (x,y,u,v)$ is torsion
free and has depth $1$. Also, $R$ is an isolated singularity. However, $\Tor_1^R(M,N)\neq 0$. This
shows that the condition (3) of \ref{vanishingiso} is critical.
\end{eg}

\end{document}